\documentclass{amsart}

\title{Multiple covers and the integrality conjecture for rational curves in
Calabi-Yau threefolds}
\author{Jim Bryan\\
Sheldon Katz\\ 
Naichung Conan Leung}
\date{\today}
\address{
Department of Mathematics\\
Tulane University\\
New Orleans, LA 70118}
\address{
Department of Mathematics\\
Oklahoma State University\\
Stillwater, OK 74078}
\address{School of Mathematics\\
University of Minnesota\\
Minneapolis, MN 55455}

\usepackage{diagrams}
\usepackage{eepic,epic}

\usepackage{amsmath,amsthm,amsfonts}
\newtheorem{thm}{Theorem}[section]

\newtheorem{lem}[thm]{Lemma}
\newtheorem{prop}[thm]{Proposition}
\newtheorem{defn}[thm]{Definition}
\newtheorem{conj}[thm]{Conjecture}

\newtheorem{Rem}[thm]{Remark}
\newenvironment{rem}{\begin{Rem}\em}{\end{Rem}}

\newcommand{\lemref}[1]{Lemma~\ref{#1}}

\newcommand{\cnums} {{\mathbf C}}          % complex numbers
\newcommand{\znums} {{\mathbf Z}}		% integers
\newcommand{\ie}{{\em i.e. }}

\newcommand{\til}[1]{\Tilde{#1}}

\renewcommand{\P}{\mathbf{P}}

\newcommand{\Def}{\operatorname{Def}}
\newcommand{\ord}[1]{\operatorname{ord}_{C_{e_{#1}}}}

\begin{document}

\begin{abstract}
We study the contribution of multiple covers of an irreducible rational
curve $C$ in a Calabi-Yau threefold $Y$ to the genus 0 Gromov-Witten
invariants in the following cases.
\begin{enumerate}
\item 
If the curve $C$ has one node and
satisfies a certain genericity condition, we prove that the contribution of
multiple covers of degree $d$ is given by
$$
\sum_{n|d}\frac{1}{n^{3}}.
$$
\item For a smoothly embedded contractable curve $C\subset Y$ we define
schemes $C_{i}$ for $1\leq i\leq l$ where $C_{i}$ is supported on $C$
and has multiplicity $i$, the number $l\in
\{1,\dots ,6 \}$ being Koll\'ar's invariant ``length''.  We prove that
the contribution of multiple covers of $C$ of degree $d$ is given by
\[
\sum _{n|d} \frac{k_{d/n}}{n^{3}}
\]
where $k_{i}$ is the multiplicity of $C_{i}$ in its Hilbert scheme
(and $k_{i}=0$ if $i>l$).
\end{enumerate}
In the latter case we also get a formula for arbitrary genus (Theorem
\ref{thm: multiple cover formula for embedded curves}).

These results show that the curve $C$ contributes an integer amount to the
so-called instanton numbers that are defined recursively in terms of the
Gromov-Witten invariants and are conjectured to be integers.
\end{abstract}

\thanks{The first author is supported by an Alfred P. Sloan Research
Fellowship and NSF grant DMS-9802612, the second author is supported by NSA
grant MDA904-98-1-0009, and the third author is supported by NSF grant
DMS-9803616.}  \maketitle \markboth{Multiple covers and
integrality}{Multiple covers and integrality}
\renewcommand{\sectionmark}[1]{}

%\tableofcontents
%\pagebreak

\section{Motivation and Results}
 From string theory and M-theory, physicists insist
on the existence of {\em instanton numbers}.  Let $Y$ be a Calabi-Yau
threefold, and $\beta\in H_2(Y,{\mathbf Z})$.  Then there is supposed to
be an integer invariant $n_{\beta}$, an instanton number, so that the
genus 0 Gromov-Witten invariant is given by
\begin{equation}
\label{multcover}
N_{\beta}=\sum_{n|\beta}\frac{n_{\beta/n}}{n^3}.
\end{equation}
A slightly more precise version of this conjecture is stated as Conjecture
7.4.5 in \cite{coxkatz}. Since Equation (\ref{multcover}) determines the
$n_{\beta }$'s recursively in terms of the $N_{\beta }$'s, one can regard the
equation as the definition of $n_{\beta }$. The integrality of $n_{\beta }$
has been verified empirically for $Y$ where the Gromov-Witten invariants
are known. However, there is no $Y$ for which the integrality of
$n_{\beta }$ has been proven to hold for all $\beta $.  M-theory suggests
an approach to an intrinsic definition of the $n_\beta$ in terms of sheaves
\cite{Go-Va}, but this has not yet been made precise mathematically.

For smoothly embedded ${\P}^1$ with normal bundle ${\mathcal O}(-1)\oplus
{\mathcal O}(-1)$, it has been proven by Aspinwall-Morrison, Voisin,
Kontsevich, Manin, Pandharipande, and Lian-Liu-Yau (see \cite{coxkatz}
\cite{Aspinwall-Morrison} \cite{man} \cite{Kont} \cite{L-L-Yau}
\cite{Voisin})\footnote{Aspinwall and Morrison's calculation predates the
definition of Gromov-Witten invariants and so their argument must be
regarded as incomplete in the present context.  It would be interesting to
find a direct comparison between their calculation and the calculations in
Gromov-Witten theory. We should also warn the reader that the theorem in
Voisin \cite{Voisin} computes a contribution of $d^{-3}$ for multiple
covers of \emph{immersed} $\P ^{1}$'s. As our results show, the
contribution of the image curve to the Gromov-Witten invariant need not be
$d^{-3}$ if the $\P^1$ is not embedded.  Voisin has only considered the
contribution of multiple covers which factor through the normalization, so
her computation does not compute the full contribution to the Gromov-Witten
invariant, as asserted in her abstract.}  that the contribution to the
genus 0 Gromov-Witten invariant of degree $d$ multiple covers is
$d^{-3}$. Thus if we imagine that all the rational curves in $Y$ are
smoothly embedded with normal bundle $\mathcal{O} (-1)\oplus \mathcal{O}
(-1)$, then $n_{\beta } $ is precisely the number of rational curves in the
class $\beta $. However this assumption is too optimistic. For example,
even if we assume Clemens' conjecture, a generic quintic 3-fold always has
nodal rational curves in degree 5 by Vainsencher \cite{Va}. Therefore, to
understand the physicists' instanton numbers and/or to extract more precise
enumerative information from the Gromov-Witten invariants, we need to
understand how non-generic curves in $Y$ contribute to the invariants.

This problem is open in almost any reasonable situation other than
$(-1,-1)$ curves.  For example: isolated smooth rational curves with a
non-generic normal bundles, and generic nodal curves. In this paper we
treat the case of any contractable smooth rational curve (which can
have normal bundle $(-1,-1),\ (0,-2)$, or $(1,-3)$) and the case of a
generic irreducible rational curve with one node.

To state our results we make the following definitions. 

A $(-1,-1)$ curve has no infinitesimal deformations so the curve is rigid;
in fact, no multiple of the curve has infinitesimal deformations so we say
the curve is super-rigid. For nodal curves, this is the notion of rigidity
that we need. It is an adaptation of a concept due to Pandharipande:

\begin{defn}[c.f. Pandharipande \cite{Pa}]\label{defn: super-rigid} A curve
$C$ in a projective variety $Y$ is called \emph{$g$-super-rigid} if for
every map $f:C'\to Y$ with $C'$ a projective curve of arithmetic genus $g$
and $f$ a local immersion, $f (C')=C$ has no infinitesimal deformations as
a stable map. A curve that is $g$-super-rigid for all $g$ is simply called \emph{super-rigid}.
\end{defn}

Super-rigidity is a genericity condition for nodal curves in the following
sense.  One can weaken the super-rigidity condition by only requiring that
the defining property holds for maps $f$ with degree at most $N$ onto $C$. Note
that for each fixed $N$, this is an open condition in any family of pairs
$(Y,C)$ and so super-rigidity is an intersection of a countable number of
open conditions. In practice, we will only need  that the condition
holds for fixed $N$. Our result for nodal curves is the following.
\begin{thm}\label{thm: formula for 1-nodal curves}
Let $C\subset Y$ be an irreducible 0-super-rigid rational curve with exactly
one node in a 3-fold $Y$ with $K_{Y}$ trivial in a neighborhood of
$C$. Then the contribution of degree $d$ multiple covers of $C$ to the
genus 0 Gromov-Witten invariant is
\[
\sum _{n|d}\frac{1}{n^{3}}.
\]
\end{thm}

Unlike the case of a smooth curve, the moduli space of stable maps that
multiply-cover a nodal curve has a number of different path components
arising from the possible ``jumping'' behavior of the map at the node. The
basic phenomenon is illustrated below for degree two maps. The moduli space
has different components depending on whether the map factors
through the normalization or not. In the figure below, points labeled by
$A$ and $B$ are mapped to corresponding points labeled $A$ and $B$ and the
normalization map identifies $A$ and $B$ to the nodal point. Notice that
the bottom map cannot be factored through the normalization (the dashed map
doesn't exist) even though the map to the nodal curve is well defined.

\vskip .25in
\setlength{\unitlength}{0.00083333in}
\begingroup\makeatletter\ifx\SetFigFont\undefined
% extract first six characters in \fmtname
\def\x#1#2#3#4#5#6#7\relax{\def\x{#1#2#3#4#5#6}}%
\expandafter\x\fmtname xxxxxx\relax \def\y{splain}%
\ifx\x\y   % LaTeX or SliTeX?
\gdef\SetFigFont#1#2#3{%
  \ifnum #1<17\tiny\else \ifnum #1<20\small\else
  \ifnum #1<24\normalsize\else \ifnum #1<29\large\else
  \ifnum #1<34\Large\else \ifnum #1<41\LARGE\else
     \huge\fi\fi\fi\fi\fi\fi
  \csname #3\endcsname}%
\else
\gdef\SetFigFont#1#2#3{\begingroup
  \count@#1\relax \ifnum 25<\count@\count@25\fi
  \def\x{\endgroup\@setsize\SetFigFont{#2pt}}%
  \expandafter\x
    \csname \romannumeral\the\count@ pt\expandafter\endcsname
    \csname @\romannumeral\the\count@ pt\endcsname
  \csname #3\endcsname}%
\fi
\fi\endgroup
{\renewcommand{\dashlinestretch}{30}
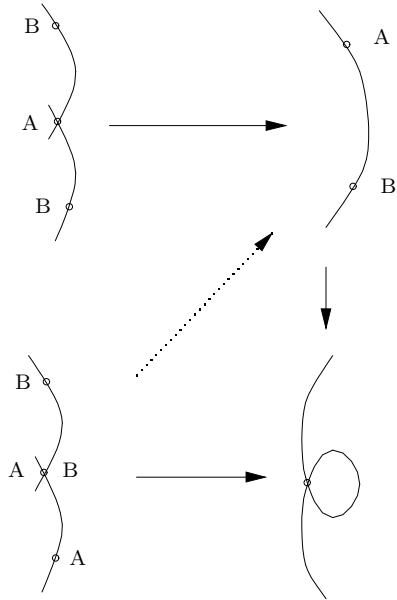
\begin{picture}(2381,3745)(-1500,-10)
\put(2151,2582){\ellipse{44}{44}}
\put(2110,3468){\ellipse{44}{44}}
\put(296,263){\ellipse{44}{44}}
\put(380,2456){\ellipse{44}{44}}
\put(222,797){\ellipse{44}{44}}
\put(307,2988){\ellipse{44}{44}}
\put(1864,734){\ellipse{44}{44}}
\put(237,1364){\ellipse{44}{44}}
\put(296,3586){\ellipse{44}{44}}
\path(802,768)(1604,768)
\blacken\path(1484.000,738.000)(1604.000,768.000)(1484.000,798.000)(1484.000,738.000)
\path(1983,2078)(1983,1697)
\blacken\path(1953.000,1817.000)(1983.000,1697.000)(2013.000,1817.000)(1953.000,1817.000)
\path(633,2963)(1730,2963)
\blacken\path(1610.000,2933.000)(1730.000,2963.000)(1610.000,2993.000)(1610.000,2933.000)
\dottedline{45}(802,1402)(1645,2288)
\blacken\path(1584.017,2180.385)(1645.000,2288.000)(1540.549,2221.743)(1584.017,2180.385)
\path(2025,1526)	(1986.685,1470.167)
	(1954.276,1421.497)
	(1927.275,1379.151)
	(1905.186,1342.293)
	(1873.754,1281.692)
	(1856.000,1233.000)

\path(1856,1233)	(1846.261,1185.646)
	(1838.990,1129.091)
	(1834.344,1066.398)
	(1832.480,1000.628)
	(1833.556,934.841)
	(1837.730,872.100)
	(1845.159,815.466)
	(1856.000,768.000)

\path(1856,768)	(1877.963,698.075)
	(1893.276,656.749)
	(1911.872,615.464)
	(1960.096,545.870)
	(2025.000,515.000)

\path(2025,515)	(2088.288,532.912)
	(2142.097,583.393)
	(2179.358,652.177)
	(2193.000,725.000)

\path(2193,725)	(2179.358,798.618)
	(2142.097,867.768)
	(2088.288,918.283)
	(2025.000,936.000)

\path(2025,936)	(1960.096,905.338)
	(1911.872,835.814)
	(1893.276,794.511)
	(1877.963,753.133)
	(1856.000,683.000)

\path(1856,683)	(1845.103,635.678)
	(1837.537,579.184)
	(1833.195,516.578)
	(1831.966,450.921)
	(1833.742,385.274)
	(1838.412,322.698)
	(1845.868,266.253)
	(1856.000,219.000)

\path(1856,219)	(1869.413,183.583)
	(1893.044,139.904)
	(1929.903,83.023)
	(1954.233,48.088)
	(1983.000,8.000)

\path(1941,3679)	(1997.237,3607.436)
	(2044.950,3544.988)
	(2084.866,3490.554)
	(2117.715,3443.030)
	(2144.224,3401.315)
	(2165.123,3364.307)
	(2193.000,3300.000)

\path(2193,3300)	(2209.018,3245.477)
	(2223.109,3179.766)
	(2234.716,3106.464)
	(2239.413,3068.091)
	(2243.281,3029.169)
	(2246.250,2990.148)
	(2248.249,2951.477)
	(2249.063,2876.987)
	(2245.165,2809.296)
	(2236.000,2752.000)

\path(2236,2752)	(2210.983,2679.122)
	(2190.767,2637.550)
	(2164.426,2590.876)
	(2131.220,2537.843)
	(2090.407,2477.192)
	(2041.247,2407.664)
	(2013.306,2369.178)
	(1983.000,2328.000)

\path(127,1526)	(154.648,1486.548)
	(180.008,1449.563)
	(224.168,1382.390)
	(260.101,1323.268)
	(288.424,1270.985)
	(309.756,1224.329)
	(324.717,1182.087)
	(333.926,1143.049)
	(338.000,1106.000)

\path(338,1106)	(331.339,1031.411)
	(320.354,989.444)
	(303.604,942.699)
	(280.590,889.929)
	(250.817,829.887)
	(213.786,761.326)
	(192.393,723.462)
	(169.000,683.000)

\path(169,895)	(191.017,853.489)
	(211.216,814.715)
	(246.412,744.755)
	(275.085,683.874)
	(297.731,630.824)
	(314.849,584.358)
	(326.935,543.231)
	(338.000,472.000)

\path(338,472)	(333.433,402.799)
	(325.326,362.246)
	(312.853,316.092)
	(295.647,263.115)
	(273.346,202.090)
	(245.585,131.793)
	(229.543,92.785)
	(212.000,51.000)

\path(212,3722)	(239.538,3682.361)
	(264.799,3645.212)
	(308.800,3577.763)
	(344.608,3518.417)
	(372.830,3465.940)
	(394.074,3419.095)
	(408.945,3376.648)
	(418.052,3337.361)
	(422.000,3300.000)

\path(422,3300)	(415.564,3225.957)
	(404.670,3184.160)
	(388.010,3137.532)
	(365.099,3084.838)
	(335.451,3024.842)
	(298.580,2956.308)
	(277.284,2918.453)
	(254.000,2878.000)

\path(254,3090)	(275.908,3048.478)
	(296.011,3009.693)
	(331.047,2939.705)
	(359.594,2878.779)
	(382.138,2825.654)
	(399.164,2779.072)
	(411.160,2737.775)
	(422.000,2666.000)

\path(422,2666)	(417.642,2597.317)
	(409.589,2556.861)
	(397.134,2510.680)
	(379.912,2457.539)
	(357.559,2396.202)
	(329.710,2325.435)
	(313.610,2286.128)
	(296.000,2244.000)

\put(2278,3468){\makebox(0,0)[lb]{\smash{{{\SetFigFont{8}{9.6}{rm}A}}}}}
\put(2321,2540){\makebox(0,0)[lb]{\smash{{{\SetFigFont{8}{9.6}{rm}B}}}}}
\put(43,1318){\makebox(0,0)[lb]{\smash{{{\SetFigFont{8}{9.6}{rm}B}}}}}
\put(338,768){\makebox(0,0)[lb]{\smash{{{\SetFigFont{8}{9.6}{rm}B}}}}}
\put(380,219){\makebox(0,0)[lb]{\smash{{{\SetFigFont{8}{9.6}{rm}A}}}}}
\put(169,2413){\makebox(0,0)[lb]{\smash{{{\SetFigFont{8}{9.6}{rm}B}}}}}
\put(98,3540){\makebox(0,0)[lb]{\smash{{{\SetFigFont{8}{9.6}{rm}B}}}}}
\put(84,2933){\makebox(0,0)[lb]{\smash{{{\SetFigFont{8}{9.6}{rm}A}}}}}
\put(0,768){\makebox(0,0)[lb]{\smash{{{\SetFigFont{8}{9.6}{rm}A}}}}}
\end{picture}
}

\vskip .25in

Thus we see that the moduli space of genus 0 stable maps that double cover
$C$ consists of two connected components (of different dimensions)
depending on whether the map factors through the normalization or not. See
Proposition \ref{prop: components of M(C0,d)} for the case of arbitrary
degree.

An isolated smoothly embedded rational curve $C$ in a Calabi-Yau
3-fold $Y$ can have a normal bundle other than $\mathcal{O} (-1)\oplus
\mathcal{O} (-1)$. In \cite{Pagota}, Reid proves that if the normal
bundle is $\mathcal{O}\oplus \mathcal{O} (-2)$, then the curve is
isolated if and only if the curve contracts; that is, there exists a
birational morphism $f:Y\to X$ with $f (C)=p\in X$ that is an
isomorphism of $Y\backslash C$ with $X\backslash p$. If a rational
curve contracts, it must be isolated, the normal bundle must be
$\mathcal{O} (-1)\oplus \mathcal{O} (-1)$, $\mathcal{O}\oplus
\mathcal{O} (-2)$, or $\mathcal{O} (1)\oplus \mathcal{O} (-3)$ and $p$
is a compound DuVal (cDV) singularity\footnote{One might have hoped that all
isolated curves contract, but this is false: Clemens has an example of
an isolated $(2,-4)$ curve where no multiple of the curve deforms
\cite{Cl} and Jim\'enez has an example of a non-contractable, isolated
$(1,-3)$-curve \cite{Ji}.}  \cite{Pi}. The singularity types of the
generic hyperplane section through $p$ were classified by
Katz-Morrison \cite{Ka-Mo} who found that they are determined by
Koll\'ar's invariant ``length''.

\begin{defn}[Koll\'ar \cite{CKM} page 95]\label{defn: length of C} Let
$C\subset Y$ be a smooth rational curve in a 3-fold $Y$ with $K_{Y}$
trivial in a neighborhood of $C$. Suppose that there exists a
contraction $\pi :Y\to X$ with $\pi (C)=p\in X$ and $\pi $ is an
isomorphism of $Y\backslash C$ onto $X\backslash p$. Define the
\emph{length} of $p$ to be the length (at the generic point of $C$) of
the scheme supported on $C$ with structure sheaf
$\mathcal{O}_{Y}/\pi ^{-1} (\mathfrak{m}_{X,p})$ where
$\mathfrak{m}_{X,p}$ is the maximal ideal of the point $p$.
\end{defn}

The length of $p$ is a number $l\in \{1,\dots ,6 \}$ and if $l=1$ then
$C$ has normal bundle $N_{C/Y}$ isomorphic to $\mathcal{O} (-1)\oplus
\mathcal{O} (-1)$ or $\mathcal{O}\oplus \mathcal{O} (-2)$ and if $l>1$
then $N_{C/Y}\cong \mathcal{O} (1)\oplus \mathcal{O} (-3) $. When $l>1$
we can define a sequence of higher order neighborhoods of $C$:

\begin{defn}\label{defn: higher order nghds of C}
Let $C\subset Y$ be as in Definition \ref{defn: length of C} and let $l$ be
the length of $p$. Let $X_{0}\subset X $ be a generic hyperplane section
through $p$ and let $Y_{0}$ be the proper transform of $X_{0}$. For
$i=1,\dots ,l$ define $C_{i}$ to be the subscheme of $Y_0$ defined by the
symbolic power $\mathcal{I}^{(i)}$ of the ideal sheaf $\mathcal{I}$
defining $C$ in $Y_0$.
\end{defn}

Note that $C_{1}=C$ and that $C_{l}$ is the scheme used in the definition
of length. The crucial property of $C_{i}$ is that under a generic
deformation of $Y$, $C_{i}$ deforms to $k_{i}$ smoothly embedded $(-1,-1)$
curves in the homology class $i[C]$ where $k_{i}$ is the multiplicity of
$C_{i}$ in its Hilbert scheme (see Proposition \ref{prop: discriminant in
PRes} and Lemma \ref{lem: intersection number equals hilbert scheme
multiplicity}).

Our main result for embedded smooth curves is the following.
\begin{thm}\label{thm: multiple cover formula for embedded curves}
Let $C\subset Y$ be a smoothly embedded contractable rational curve in a
3-fold $Y$ with $K_{Y}$ trivial in a neighborhood of $C$.  Then the
contribution of degree $d$ multiple covers of $C$ to the genus $g$
Gromov-Witten invariants is
\[
\sum _{n|d} k_{d/n}\frac{|B_{2g}|n^{2g-3}}{2g\cdot(2g-2)!}
\]
where $B_{2g}$ is the $2g$-th Bernoulli number
(c.f. \cite{Fa-Pa} Theorem 3) and $k_{i}$ is the multiplicity of $C_{i}$
in its Hilbert scheme if $i\leq l$ and $k_{i}=0$ otherwise. In particular,
the genus 0 formula is
\[
\sum _{n|d} \frac{k_{d/n}}{n^{3}}.
\]
\end{thm}
For example, if $l=1$ then $C$ is a $(-1,-1)$ or $(0,-2)$ curve and
$k_{1}$ coincides with Reid's invariant ``width'' (Definition 5.3 of
\cite{Pagota}).

Our technique for proving Theorem \ref{thm: formula for 1-nodal
curves} is to relate the invariants of a super-rigid 1-nodal curve to
the invariants of a contractable chain of rational curves via a method
similar to that used by Bryan-Leung to study contributions of multiple
covers of nodal curves in $K3$ surfaces (Section 5 of
\cite{Br-Le1}). The invariants of a contractable chain of curves (as
well as the contractable curve in Theorem \ref{thm: multiple cover
formula for embedded curves}) are then computed by deforming a
neighborhood of the curve. The deformation is constructed and analyzed
by employing the beautiful subject of deformations of DuVal
singularities and their simultaneous (partial) resolutions. The
subject goes back to Brieskorn \cite{Br1}\cite{Br2}\cite{Br3} and has
been used to study 3-fold singularities with small resolutions by
Pinkham \cite{Pi} and others (c.f. \cite{Fr}\cite{Ka-Mo}\cite{Pagota}).

It is crucial that the deformations we construct in this way can be blown
down to an open subset (in the analytic topology) of an affine
variety. This enables us to conclude that all the complete curves in the
deformed 3-fold lie in the exceptional set of the blowdown.

Unfortunately, we do not know how to generalize these techniques to curves
with more than one node; in fact, we do not even have a reasonable
conjectural formula (see Remark \ref{rem: problem with more nodes}). There
are also problems generalizing the proof of Theorem \ref{thm: formula for
1-nodal curves} to higher genus (see Remark \ref{rem: problem with higher
genus}). However, considerations of M-theory \cite{Go-Va} lead us to
the following conjecture:
\begin{conj}\label{conj: higher genus for 1-nodal curves}
Let $C\subset Y$ be an irreducible super-rigid rational curve with exactly
one node in a 3-fold $Y$ with $K_{Y}$ trivial in a neighborhood of
$C$. Then the contribution of degree $d$ multiple covers of $C$ to the
genus $g$ Gromov-Witten invariants is
\[
\sum _{n|d}\left(\frac{|B_{2g}|n^{2g-3}}{2g (2g-2)!}+\frac{\delta
_{1,g}}{n} \right)
\]
where $\delta _{1,g}=1$ if $g=1$ and $\delta _{1,g}=0$ if $g\neq 1$.
\end{conj}

Via M-theory, genus $g$ invariants $n^g_\beta$ can be assigned for all
$g\ge0$ and all homology classes $\beta$, from which the 0-point
Gromov-Witten invariants can be conjecturally computed by an explicit
formula generalizing equation (\ref{multcover}), with $n^0_\beta=n_\beta$
(see \cite{Go-Va}).  Let us consider these invariants (if they exist) to be
{\em defined\/} by this formula (c.f. \cite{Pa}). We can then discuss the
contributions of a curve $C$ to $n^{g}_{d[C]}$.  Theorem~\ref{thm: formula
for 1-nodal curves} proves that the contribution of $C$ to $n_{d[C]}$ is
$1$ for all $d\ge 1$, and Conjecture~\ref{conj: higher genus for 1-nodal
curves} asserts that the contribution to $n^1_{d[C]}$ is also $1$
for all $d\ge 1$ while the contribution to $n^g_{d[C]}$ is $0$ for $g\ge2$.

In this language, \cite{Go-Va} also gives a conjectural formula for these
invariants in terms of moduli spaces of $U(d)$ bundles on $C$.  These
moduli spaces have been computed in \cite{Teodorescu-thesis}, and we have
checked that this computation leads to $C$ contributing $1$ to $n_{d[C]}$
for all $d$, as predicted.

Our paper is organized as follows. In Section \ref{sec: nghd of contractable
curve and defs} we study the neighborhood of a contractable curve in a
Calabi-Yau 3-fold. We explicitly construct deformations of the neighborhood
in order to compute its contribution to the Gromov-Witten invariants. In
Section \ref{sec: relating nodes to chains}, we relate the Gromov-Witten
contribution of a 0-super-rigid nodal curve to the contribution of
contractable chains of curves.  As an example of Theorem \ref{thm: multiple
cover formula for embedded curves}, we work out the length $2$ case in
detail in Section \ref{sec: apps and examples}, computing $k_{1}$ and
$k_{2}$ directly from explicit equations for the 3-fold.

\section{The neighborhood of a contractable curve and its 
deformations}\label{sec: nghd of contractable curve and defs}

In this section we use the deformation invariance property of Gromov-Witten
invariants to compute them for a neighborhood of a contractable curve in a
local Calabi-Yau 3-fold. This is achieved by deforming the neighborhood so
that all the curves in the neighborhood are smooth $(-1,-1)$ curves
(c.f. Friedman \cite{Fr} pg. 678--679). We focus on the two cases of
interest: a generic linear chain of curves that contracts to a $cA_{n}$
singularity and a single contractable smooth curve which necessarily must
contract to a $cA_{1}$, $cD_{4}$, $cE_{6}$, $cE_{7}$, or $cE_{8}$
singularity. We will relate the former case to the invariants of a nodal
curve in Section \ref{sec: relating nodes to chains}.

We begin by outlining the theory of deformations of DuVal
singularities and their simultaneous resolutions. This is due to
Brieskorn \cite{Br1}\cite{Br2}\cite{Br3} and, in a refined form that
we need, Katz-Morrison (see Section 3, especially Theorem 1 of
\cite{Ka-Mo}). The application of these ideas to Gorenstein 3-fold
singularities with small resolutions was pioneered by Pinkham
\cite{Pi} and Reid \cite{Pagota}.

Let $C\subset Y$ be a rational curve (not necessarily irreducible) in a
3-fold $Y$ with $K_{Y}$ trivial in a neighborhood of $C$ and $\pi :Y\to X$
is a birational morphism such that $\pi (C)=p\in X$ and $\pi |_{Y\backslash
C}$ is an isomorphism onto $X\backslash p$. We consider an analytic
neighborhood of $p$ (still denoted $X$) and its inverse image under $\pi $
(still denoted $Y$). By a lemma of Reid \cite{Pagota} (1.1,1.14), the
generic hyperplane section through $p$ is a surface $X_{0}$ with an
isolated rational double point, and the proper transform of $X_{0}$ is a
partial resolution $Y_{0}\to X_{0}$ (\ie the minimal resolution $Z_{0}\to
X_{0}$ factors through $Y_{0}\to X_{0}$).

The partial resolution $Y_{0}\to X_{0}$ determines combinatorial data
$\Gamma _{0}\subset \Gamma $ consisting of an ADE Dynkin diagram
$\Gamma $ (the type of the singularity $p$) and a subgraph $\Gamma
_{0}$ (the dual graph of the exceptional set of $Y_{0}$).

Let $\mathcal{Z}\to \Def (Z_{0})$, $\mathcal{Y}\to \Def (Y_{0})$, and
$\mathcal{X}\to \Def (X_{0})$ be semi-universal deformations of
$Z_{0}$, $Y_{0}$, and $X_{0}$. Following \cite{Ka-Mo}, there are
identifications
\begin{align*}
\Def (Z_{0})&\cong U=:\operatorname{Res} (\Gamma )\\
\Def (Y_{0})&\cong U/W_{0}=:\operatorname{PRes} (\Gamma ,\Gamma _{0})\\
\Def (X_{0})&\cong U/W=:\Def (\Gamma )
\end{align*}
where $U$ is the complex root space associated to $\Gamma $ and $W$ is
its Weyl group. $W_{0}\subset W$ is the subgroup generated by
reflections of the simple roots corresponding to $\Gamma -\Gamma
_{0}$. Deformations of $Z_{0}$ or $Y_{0} $ can be blown down to give
deformations of $X_{0}$ (\cite{Wahl} Theorem 1.4) and the induced
classifying maps are given by the natural maps $U\to U/W$ and
$U/W_{0}\to U/W$ under the above identifications.

Via the defining equation for the hyperplane section $X_{0}$, we can
view $X$ as the total space of a 1-parameter family $X_{t}$ defined by
a classifying map
\[
g:\Delta \to \Def (\Gamma ).
\]     
Similarly, we get a compatible family $Y_{t}$ given by a map
\[
f:\Delta \to \operatorname{PRes} (\Gamma ,\Gamma _{0}).
\]
That is, we get the diagram
\[
\begin{diagram}
\mathcal{Z}&\rTo^{\tilde{\sigma}}&\mathcal{Y}&\rTo^{\tilde{\rho}}&\mathcal{X}\\
\dTo&&\dTo&&\dTo\\
\operatorname{Res} (\Gamma )&\rTo^{\sigma }&\operatorname{PRes} (\Gamma ,\Gamma _{0})&\rTo^{\rho }&\Def (\Gamma )\\
&&\uTo>{f}&\ruTo_{g}&\\
&&\Delta &&
\end{diagram}
\]
where $Y$ is the pullback of $\mathcal{Y}$ by $f$ and $X$ is the
pullback of $\mathcal{X}$ by $g$.

Our aim is to deform $Y$ by deforming the map $f$ and to describe the
curves in the resulting space. To do this we need to understand the
discriminant loci of $\operatorname{PRes} (\Gamma ,\Gamma _{0})$ and the
curves in the corresponding surfaces. We first remark that if
$\til{f}:\Delta \to \operatorname{PRes} (\Gamma ,\Gamma _{0})$ is a generic
deformation of $f$ (so in particular $\til{f}$ is transverse to the
discriminant locus), then $\til{Y}:=\til{f}^{-1} (\mathcal{Y})$ is
smooth. This is true because the singularities of $\mathcal{Y}$ all lie
over a set of at least codimension 2 in $\operatorname{PRes} (\Gamma
,\Gamma _{0})$.  Ultimately, we show that all the curves in $\til{Y}$ are
smooth $(-1,-1)$ curves. If we then know the number of such curves and
their homology class, then the Gromov-Witten invariants of $\til{Y}$ can be
easily computed from the following multiple cover formula of Faber and
Pandharipande.

\begin{lem}[Faber-Pandharipande \cite{Fa-Pa}]\label{lem: multiple cover formula in genus g for (-1,-1) curves} 
If $C$ is a smoothly embedded $(-1,-1)$ curve in a Calabi-Yau 3-fold, then
the contribution of degree $d$ multiple covers of $C$ to the genus $g$
Gromov-Witten invariants is given by
\[
\frac{|B_{2g}|d^{2g-3}}{2g\cdot (2g-2)!}.
\]
\end{lem}

The exceptional curve of the minimal resolution $Z_{0}\to X_{0}$ is a
rational curve whose dual graph is $\Gamma $. Thus the components of the
exceptional divisor are rational curves $C_{e_{i}}$ which are in one to one
correspondence with the simple roots $e_{1},\dots ,e_{n}$ ($n$ is the rank
of the root system). A generic deformation of $Z_{0}$ has no complete
curves, but the curves $C_{e_{i}} $ and various combinations of $C_{e_{i}}$
will deform in certain codimension one families of deformations. To make
this precise, we define the \emph{discriminant locus} $\mathfrak{D} \subset
\operatorname{Res} (\Gamma )$ to be those $t\in \operatorname{Res} (\Gamma
)$ such that the corresponding surface $\mathcal{Z}_{t}$ has a complete
curve. The following proposition is essentially part 3 of Theorem 1 in
\cite{Ka-Mo}.
\begin{prop}\label{prop: components of disc locus of Res}
The irreducible components of the discriminant divisor $\mathfrak{D} \subset
\operatorname{Res} (\Gamma )$ are in one to one correspondence with the
positive roots of $\Gamma $. Under the identification of
$\operatorname{Res} (\Gamma )$ with the complex root space $U$, the
component $\mathfrak{D} _{v}$ corresponding to the positive root $v=\sum
_{i=1}^{n}a_{i}e_{i}$ is $v^{\perp
}\subset U $, \ie  the hyperplane perpendicular to $v$.

Moreover, $\mathfrak{D} _{v}$ corresponds exactly to those deformations of
$Z_{0}$ in $\mathcal{Z}$ to which the curve 
\[
C_{v}:=\bigcup
_{i=1}^{n}a_{i}C_{e_{i}}
\]
lifts. For a generic point $t\in \mathfrak{D} _{v}$,
the corresponding surface $\mathcal{Z}_{t}$ has a single smooth $-2$ curve
in the class $\sum _{i=1}^{n}a_{i}[C_{e_{i}}]$ thus there is a small
neighborhood $B$ of $t$ such that the restriction of $\mathcal{Z}$ to $B$
is isomorphic to a product of $\cnums ^{n-1}$ with the semi-universal
family over $\operatorname{Res} (A_{1})$.
\end{prop}

Katz and Morrison prove this (among other things) by constructing explicit
equations for $\mathcal{Z}$. Although we do not need the explicit form of
these equations for the proofs of our theorems, they are useful for
computing the values of $k_{i}$ in concrete cases (see Section \ref{sec:
apps and examples}).

Since the Weyl group acts transitively on the positive roots, we see that
the discriminant locus in $\operatorname{Def} (\Gamma )$, which is the
image of $\mathfrak{D} $, is irreducible. It corresponds to those
deformations of $X_{0}$ which are singular.

The two cases relevant to this paper are when $\Gamma =\Gamma _{0}=A_{n}$
(so $\operatorname{PRes} (\Gamma ,\Gamma _{0})=\operatorname{Res} (\Gamma
)$) and when $\Gamma _{0}$ is a single vertex. In this latter case, $\Gamma
$ must be one of $A_{1}$, $D_{4}$, $E_{6}$, $E_{7}$, or $E_{8}$ and $\Gamma
_{0}$ must be the central vertex or, in the case of $E_{8}$, it can also be
the vertex one away from the center on the long branch. These six
possibilities correspond to the six possible values of the length (this is
the main theorem of \cite{Ka-Mo}).

\subsection{The case of a contractable smooth $\P ^{1}$.}  Assume then that
$\Gamma _{0}$ is a single vertex so that $(\Gamma ,\Gamma _{0})$ is one of
the above six possibilities. Let $D\subset \operatorname{PRes} (\Gamma
,\Gamma _{0})$ denote the discriminant locus, \ie the image of
$\mathfrak{D} $ under the quotient of $W_{0}$. The corresponding surfaces
are those which are singular or contain complete curves (or both). We write
$D=D^{sing}\cup D^{curv}$ where the surfaces corresponding to the points in
$D^{curv}$ contain a complete curve and $D^{sing}$ is the union of the
remaining components. A generic point of $D^{curv}$ corresponds to a smooth
surface with a smooth (-1,-1) curve and a generic point of $D^{sing}$
corresponds to a surface with no curves and a single $A_{1}$ singularity.

For the purposes of computing the Gromov-Witten invariants of $Y$ we need
to understand the components of $D^{curv}$. Since the family $\mathcal{Y}$ is
obtained from $\mathcal{Z}$ by simultaneously blowing down the curves in
$\mathcal{Z}_{t}$ corresponding to $\Gamma -\Gamma _{0}$ followed by
taking the quotient by $W_{0}$, we can study $D^{curv}$ by studying the orbit
structure of $W_{0}$ acting on $\mathfrak{D} =\cup \mathfrak{D} _{v}$.

\begin{prop}\label{prop: discriminant in PRes}
The divisor $D^{curv}\subset \operatorname{PRes} (\Gamma ,\Gamma _{0})$
consists of $l$ irreducible components $D_{1},\dots ,D_{l}$.  Moreover, if
$t$ is a generic point in $D_{i}$ then the surface $\mathcal{Y}_{t}$ has a
smooth $-2$ curve in the class $i[C]$.  These curves form the fibers of
a divisor $\mathcal{C}_i$ inside the restriction  of $\mathcal{Y}$ over
$D_i$.  The scheme-theoretic fiber of $\mathcal{C}_{i}$ over the origin is
precisely the subscheme $C_i\subset Y_0$.
\end{prop}

\textsc{Proof:} Consider the action of $W_{0}$ on the roots $\{v \}$ and
the induced action on the hyperplanes $\{\mathfrak{D} _{v} \}$. Let $e_{1}$
denote the simple root corresponding to $\Gamma _{0}$. The coefficient
$\alpha _{1}$ of any root $v=\sum _{i=1}^{n}\alpha _{i}e_{i}$ is preserved
by the generators of $W_{0}$ and so it is an invariant of the orbits of
$W_{0}$ acting on $\{\mathfrak{D} _{v} \}$. For orbits with $\alpha _{1}\neq 0$
it turns out that $\alpha _{1}$ is a complete invariant:

\begin{lem}\label{lem: a1 coef is a complete invariant of W0 orbits}
If $v=\sum _{i=1}^{n}\alpha _{i}e_{i}$ is a root with $\alpha _{1}\neq 0$,
then $v'=\sum _{i=1}^{n}\alpha '_{i}e_{i}$ is in the $W_{0}$-orbit of $v$
if and only if $\alpha' _{1}=\alpha _{1}$.
\end{lem}

Since the groups and root systems are finite and there are only a finite
number of cases, this can easily be checked by hand. It is a fun exercise
(really!) which we encourage the reader to carry out.

The $W_{0}$ orbits of $\{\mathfrak{D} _{v} \}$ with $\alpha _{1}\neq 0$ are
exactly the components of $D^{curv}$ and the $W_{0}$ orbits of
$\{\mathfrak{D} _{v} \}$ with $\alpha _{1}=0$ are the components of
$D^{sing}$. By inspection, the possible non-zero values of $\alpha _{1}$
are $1,\dots ,l$, and so define $D_{k}$ to be the image of $\{\mathfrak{D}
_{v}:v=\sum \alpha _{i}e_{i}\text{ has }\alpha _{1}=k \}$ in 
$\operatorname{PRes} (\Gamma,\Gamma_1)$. Since the
surface $\mathcal{Z}_{t}$ for $t\in \mathfrak{D} _{v}$ has a curve in the
class $\sum _{i=1}^{n}\alpha _{i}[C_{e_{i}}]$, the surface
$\mathcal{Y}_{t'}$ (where $t\mapsto t'$) has a curve in the class $\alpha
_{1}[C]$. For generic $t\in \mathfrak{D} _{v}$ this curve is a single
smooth $(-1,-1) $ curve and $\mathcal{Z}_{t}\cong \mathcal{Y}_{t'}$.

For the claim about $C_i$, note that after pulling back to $\operatorname{Res}
(\Gamma)$, we can pick a discriminant component $\mathfrak{D}_v$ mapping
to $D_i$, and a corresponding (possibly reducible) divisor 
$\mathcal{C}_v$ in the restriction
of $\mathcal{Z}$ to $\mathfrak{D}_v$ whose fibers are the (possibly reducible)
curves corresponding to the discriminant component $\mathfrak{D}_v$.  
The divisor $\mathcal{C}_i$ is the
scheme-theoretic image of $\mathcal{C}_v$ under the restriction
of $\tilde{\sigma}$ over the relevant
discriminant components.  The fiber $C_v$ of $\mathcal{C}_v$ over the origin
is itself an exceptional $\mathrm{ADE}$ configuration, and its
scheme structure is defined by the pullback of the maximal ideal of the
singularity being resolved.  Writing its
cycle class as $C_v=\sum_{j\ge1}\alpha_jC_{e_j}$, we have that $\alpha_1=i$. 
We have that $C_v$ is scheme-theoretically defined by the ideal of
functions vanishing to order $\alpha_j$ at each $C_{e_{j}}$~\cite{BPV}.
We now look at the scheme-theoretic image of $C_v$ after contracting
the $C_{e_{j}}$ for $j>1$ via $\tilde{\sigma}$.  Its ideal consists of those 
functions $f$ such that $\ord{j}\tilde{\sigma}^*(f)\ge\alpha_j$, where
$\ord{j}$ denotes the order of vanishing along $C_{e_j}$.  It remains to
show that if $f\in\mathcal{I}^{(i)}$, then $\ord{j}\tilde{\sigma}^*(f)\ge
\alpha_j$ for all $j$.  To expedite this task we use the following lemma.
\begin{lem}
\label{Lemma:subadditivity} Let $f$ be a function on a neighborhood of $C$
in $Y_0$.  Put $m_k=\ord{j}\tilde{\sigma }^*(f)$. Then for each
$\tilde{\sigma }$-exceptional curve $C_{e_k}$ we have $2m_k\ge \sum_j m_j$.
\end{lem}
The lemma follows immediately from the obvious assertion: for all
$\tilde{\sigma }$-exceptional curves $D$, we have $\tilde{\sigma }^*((f))
\cdot D \le 0$.

It is now a simple matter to check in each case that the inequalities of
\lemref{Lemma:subadditivity} together with $m_1\ge\alpha_1=i$ imply that
$m_k\ge\alpha_k$ for all $k$.
\qed

\begin{Rem}
The proof of Proposition~\ref{prop: discriminant in PRes} shows that the
scheme theoretic exceptional set of $Y_0\to X_0$ is defined by the ideal
$\mathcal{I}^{(l)}$.
\end{Rem}

\smallskip
We now wish to compute the Gromov-Witten invariants of $Y$ in the class
$d[C]$.  Let $\til{f}:\Delta \to \operatorname{PRes} (\Gamma ,\Gamma _{0})$
be a small generic deformation of $f$ and define $\til{Y}$ to be the total
space of the induced one parameter family of surfaces, $\til{f}^{-1}
(\mathcal{Y})$.

\smallskip
It makes sense to speak of the Gromov-Witten invariants of a local manifold
like $Y$ since for any closed 3-fold $\overline{Y}$ containing $Y$, the
moduli space of (non-constant) stable maps has a connected component
consisting of maps whose image lies entirely in $Y$. This is because a
curve that lies in $Y$ must be in the exceptional set of $Y\to X$ since $X$
is an open set in an affine variety. We understand the ``Gromov-Witten
invariants of $Y$'' to mean the virtual fundamental cycle restricted to
that component. Alternatively, one can apply the analytic definition of the
Gromov-Witten invariants directly to $Y$ since the above argument shows
that the moduli space of stable maps is compact, the (analytic) virtual
class constructions of Li and Tian apply \cite{Li-Tian}.  

\begin{lem}\label{lem: GW of Y and tilY are well defined}
The Gromov-Witten invariants of $\til{Y}$ are well defined and equal to the
Gromov-Witten invariants of $Y$.
\end{lem}

By the above
argument, all of the complete curves of our deformations of $Y$ lie in the
exceptional set of $\mathcal{Y}\to \mathcal{X}$ which in turn lies over
$D^{curv}\subset \operatorname{PRes} (\Gamma ,\Gamma _{0})$. Let
\[
F:\Delta \times [0,\epsilon ]\to \operatorname{PRes} (\Gamma ,\Gamma _{0})
\]
be a smooth map so that $f_{s}:=F|_{\Delta \times \{s \}}$ is analytic,
$f_{0}=f$, and $F$ is transverse to $D^{curv}$. Furthermore, we may assume,
by making $\epsilon $ smaller if necessary, that $F^{-1} (D^{curv})$ is
bounded away from $\partial \Delta $. Let $\til{f}=f_{\epsilon }$. The
parameterized moduli spaces are then compact and so the deformation
invariance arguments in \cite{Li-Tian} apply.\qed

Since $\til{f}$ is transverse to $D^{curv}$, by Proposition \ref{prop:
discriminant in PRes}, we see that the complete curves of $\til{Y}$ are all
smooth $(-1,-1)$ curves in the classes $[C]$, $2[C]$,\dots , $l[C]$. The
number of irreducible curves in the class $i[C]$ is the intersection number
of $\Delta $ with $D_{i}$, \ie the cardinality of $\til{f}^{-1}
(D_{i})$. Thus Theorem \ref{thm: multiple cover formula for embedded
curves} follows easily from Lemma \ref{lem: multiple cover formula in genus
g for (-1,-1) curves} once we show that the intersection numbers of $\Delta
$ with $D_{i}$ are equal to the Hilbert scheme multiplicities.
\begin{lem}\label{lem: intersection number equals hilbert scheme multiplicity}
The intersection number $\# \{\til{f}^{-1} (D_{i}) \}$ is equal to $k_{i}$,
the multiplicity of $C_{i}$ (Definition \ref{defn: higher order nghds of
C}) in its Hilbert scheme.
\end{lem}
\textsc{Proof:} For each point in $D_{i}$, the corresponding surface has a
unique irreducible complete curve. Let $\mathcal{C}_{i}\to D_{i}$ be the
corresponding family of curves. We claim this family is a flat family of
rational curves $\mathcal{C}_{i}\subset \mathcal{Y}$.

{}From the explicit construction in \cite{Ka-Mo}, we see that
$\mathcal{C}_{i}\to D_{i}$ is projective so for flatness, it suffices to
check that the genus and degree of the fibers are constant where any
convenient polarization can be chosen for the purpose of defining the
degree.  The fibers are all rational, since they are the exceptional set of
the map of a surface to a surface with a rational singularity. To see that
the degree is constant, consider the pullback of $\mathcal{C}_{i}\to D_{i}$
to $\operatorname{Res} (\Gamma )$. We get a family $\sigma^{-1}
(\mathcal{C}_{i}) \to \mathfrak{D}_{i}$ where $\mathfrak{D}_{i}=\sigma
^{-1} (D_{i}) $ is the union of the components $\mathfrak{D}_{v}$ mapping
to $D_{i}$. It is enough to show that $\sigma^{-1} (C_{i})\to
\mathfrak{D}_{i}$ has constant fiberwise degree. By restricting
$(\tilde{\rho} \circ \tilde{\sigma} )^{-1} (\mathcal{X})$, $\tilde{\sigma}
^{-1} (\mathcal{Y})$, and $\mathcal{Z}$ to $\mathfrak{D}_{i}$, we get the
families of singular, partially resolved, and minimally resolved surfaces
over $\mathfrak{D}_{i}$ which we denote by $\mathcal{X}'$, $\mathcal{Y}'$,
and $\mathcal{Z}'$ respectively.  From \cite{Ka-Mo}, these surfaces all
have projective compactifications and hence polarizations
$L_{\mathcal{X}'}$, $L_{\mathcal{Y}'}$, and $L_{\mathcal{Z}'}$. We need to
show that $\sigma^{-1} (\mathcal{C}_{i})\subset \mathcal{Y}'$ has constant
degree, \ie the fibers of $\sigma ^{-1} (\mathcal{C}_{i})$ have constant
intersection with $L_{\mathcal{Y}'}$, which we will refer to as `` $\sigma
^{-1} (\mathcal{C}_{i})\cdot L_{\mathcal{Y}'}$ is constant on the
fibers''. Let $\mathcal{C}\subset \mathcal{Z}'$ be the exceptional set for
the blowdown $\mathcal{Z'\to X'}$ and let $\mathcal{C=C'+C''}$ where
$\mathcal{C''}$ is the exceptional set for $b:\mathcal{Z'\to Y'}$.  Then
the exceptional set $\mathcal{C}\to \mathfrak{D_{i}}$ is flat.  This was
written explicitly for all $v$ in \cite{Ka-Mo} in the cases of $A_n$ (page
469) and $D_n$ (page 473).  This also holds true for $E_n$ as can be
checked directly case by case.  It follows that $\mathcal{C}\cdot
L_{\mathcal{Z}'}$ is constant on fibers.  Furthermore, $L_\mathcal{Z}'$ can
be chosen to be the pullback of $L_\mathcal{Y}'$ plus an appropriate
multiple $m$ of $\mathcal{C}''$.  Since $b (\mathcal{C}'')$ is of lower
dimension, $b_{*} (\mathcal{C''})=0$, and we have as an identity on the
fibers
\[
\sigma ^{-1} (\mathcal{C}_{i})\cdot L_{\mathcal{Y'}}=
b_{*} (\mathcal{C'})\cdot
L_{\mathcal{Y'}}=b_{*} (\mathcal{C})\cdot L_{\mathcal{Y}'}=
\mathcal{C}\cdot \left(L_{\mathcal{Z}'}-m\mathcal{C}''\right),
\]
which is constant since the intersection number of the fibers of
$\mathcal{C}$ and $\mathcal{C}''$ is constant.
Thus $\mathcal{C}_{i}\to D_{i}$ is flat.

Now consider an analytic family $f_{\epsilon }$ of deformations of $f$
parameterized by a disk $B$, \ie an analytic map
\[
F:\Delta \times B \to \operatorname{PRes} (\Gamma ,\Gamma _{0})
\]
where $F|_{\Delta \times \epsilon }$ is denoted $f_{\epsilon }$. Assume
that $\til{f}$ is given by $f_{1}$. We get a finite map $\til{B}\to B$
whose degree is the intersection number $\# \{\til{f}^{-1} (D_{i}) \}$ by
projecting $F^{-1} (D_{i})$ to $B$. By the flatness of $\mathcal{C}_{i}\to
D_{i}$, we get a flat family of curves $F^{-1} (\mathcal{C}_{i})\to
B$. This family sit inside the family of 3-folds $F^{-1} (\mathcal{Y})\to
B$ (whose fiber over $\epsilon $ is the 3-fold $f^{-1}_{\epsilon }
(\mathcal{Y})$) and so determines a subscheme of the relative Hilbert
scheme (or Douady space) of $F^{-1} (\mathcal{Y})/B$ that is finite over
$B$.This family is isomorphic to $\til{B}\to B$ and since its degree is
exactly the multiplicity of $C_{i}$ in $Y=f^{-1}_{0} (\mathcal{Y})$, the
lemma is proved.

\subsection{The case of a generic contractable $A_{n}$
curve.} \label{subsec: generic contractable An}

We now assume that $Y$ is a neighborhood of a generic, contractable
chain of $n$ rational curves $C=C_{e_{1}}\cup \dots \cup C_{e_{n}}$. Thus
we are in the case where $\Gamma _{0}=\Gamma =A_{n}$ and $Y=f^{-1}
(\mathcal{Y})$ where $f:\Delta \to \operatorname{Res} (A_{n})$. Our
genericity assumption is that $f$ is transverse to the components $\mathfrak{D}
_{ij}\subset \operatorname{Res} (A_{n})$ of the discriminant locus. That
is, although $f$ intersects the discriminant locus only at the origin, it
meets each hyperplane $\mathfrak{D} _{ij}$ transversely. Here the hyperplanes
$\{\mathfrak{D} _{ij}:1\leq i\leq j\leq n \}$ are perpendicular to the roots
$v_{ij}=\sum _{k=i}^{j}e_{k}$. This condition is related to the
non-existence of infinitesimal deformations:
\begin{lem}\label{lem: f transverse is equiv to no inf defs}
Let $f:\Delta \to \operatorname{Res} (A_{n})$ be as above so that $f$
intersects each hyperplane $\mathfrak{D} _{ij}$ transversely exactly once at
$0\in \operatorname{Res} (A_{n})$. Let $Y=f^{-1} (\mathcal{Y})$ and let
$C_{e_{1}}\cup \dots \cup C_{e_{n}}=C\subset Y$ be the corresponding
contractable curve. For any $1\leq i\leq j\leq n$, let $C_{ij}=C_{e_{i}}\cup
\dots \cup C_{e_{j}}$ and let $N_{ij}$ be the normal bundle of $C_{ij}$ in
$Y$. Then $H^{0} (C_{ij},N_{ij})=0$ for all $i$ and $j$; \ie $C_{ij}$ has
no infinitesimal deformations. Conversely, if $H^{0} (C_{ij},N_{ij})\neq
0$, then $f$ meets $\mathfrak{D} _{ij}$ non-transversely.
\end{lem}
\textsc{Proof:} This follows directly from the description of deformations
of the curves $C_{ij}$ in Chapter~6 of \cite{Zerger-thesis},
or from the characterization of $\mathfrak{D}
_{v}$ in terms of deformations of $C_{v}$ that was given in Proposition
\ref{prop: components of disc locus of Res}.  Suppose that for some $i$ and
$j$, $C_{ij}$ has a non-trivial infinitesimal deformation. It defines a
morphism $h:\mathcal{C}\to Y$ where $\pi :\mathcal{C}\to \mathbf{I}$ is a
curve over $\mathbf{I}=\operatorname{Spec} (\cnums [\epsilon ]/ \epsilon
^{2})$ and $h|_{\pi ^{-1} (0)}:C_{ij}\hookrightarrow Y$. Since
$\mathfrak{D} _{ij}\subset \operatorname{Res} (A_{n})$ classifies those
deformations to which $C_{ij}$ lifts, the map $\mathcal{C}\to Y$ is induced
by a classifying map $\mathbf{I}\to \mathfrak{D} _{ij}$ so we get a
commutative diagram
\[
\begin{diagram}
\mathcal{C}&\rTo&Y&\rTo&\mathcal{Y}&\lInto&\mathcal{Y}|_{\mathfrak{D} _{ij}}&\lTo&\mathcal{C}\\
\dTo&&\dTo&&\dTo&&\dTo&&\dTo\\
\mathbf{I}&\rTo&\Delta &\rTo&\operatorname{Res} (A_{n})&\lInto&\mathfrak{D}
_{ij}&\lTo&\mathbf{I}.
\end{diagram}
\]
and so $\Delta $ lies in $\mathfrak{D} _{ij}$ to first order and we
conclude that $f$ is not transverse to $\mathfrak{D} _{ij}$. Conversely, if
$f$ is not transverse to $\mathfrak{D} _{ij}$, then the reverse of the
above argument gives an infinitesimal deformation of $C_{ij}$.\qed

We can now compute the Gromov-Witten invariants of $Y$ as in the previous
case by deforming $Y=f^{-1} (\mathcal{Y})$ to $\til{Y}=\til{f}^{-1}
(\mathcal{Y})$. The Gromov-Witten invariants of $Y$ and $\til{Y}$ are well
defined and equal by the same argument as Lemma \ref{lem: GW of Y and tilY
are well defined} and all the curves of $\til{Y}$ are smooth $(-1,-1)$
curves arising from the transverse intersection of $\til{f}$ with $\mathfrak{D}
_{ij}$ at generic points. Thus we see that there is exactly one $(-1,-1)$
curve in the class $\sum _{k=i}^{j}[C_{e_{k}}]$ for each $1\leq i\leq j\leq
n$ and these are all of the complete curves in $\til{Y}$. Combining this
observation with Lemma \ref{lem: multiple cover formula in genus g for
(-1,-1) curves} we get the following proposition:
\begin{prop}\label{prop: GW of chain}
Assume that $d_{i}>0$ for $i=1,\dots ,n$. Let $N_{g} (d_{1},\dots ,d_{n})$
denote the contribution of the genus $g$ Gromov-Witten invariant of $Y$ in
the class $\sum _{i=1}^{n}d_{i}[C_{e_{i}}]$. Then
\[
N_{g} (d_{1},\dots ,d_{n})=
\begin{cases}
\frac{|B_{2g}|}{2g (2g-2)!}d^{2g-3}&\text{if $d_{k}=d$ for all $k$}\\
0&\text{otherwise.}
\end{cases}
\]
\end{prop}

It is straightforward to formulate and prove a generalization of 
Proposition~\ref{prop: GW of chain} to the non-generic case, where each 
$C_{ij}$ is assumed isolated but not necessarily infinitesimally isolated.
The multiplicity of each $C_{ij}$ in its Hilbert scheme appears in the
formula.

\section{Relating nodal curves to chains.}\label{sec: relating nodes to chains}
In this section we relate the Gromov-Witten invariants of a
super-rigid nodal rational curve in a Calabi-Yau 3-fold to the
Gromov-Witten invariants of a neighborhood of a contractable chain of curves
in a Calabi-Yau 3-fold. 

Let $C_{0}\subset \overline{Y}$ be a super-rigid, irreducible, rational
curve $C_{0}$ with exactly one node in a projective 3-fold
$\overline{Y}$. We consider a neighborhood (in the analytic topology) of
$C_{0}$ which we denote by $Y$ that is sufficiently small enough so that
all non-constant stable maps to $\overline{Y}$ whose image is contained in
$Y$ have image $C_{0}$. Then the moduli space $\overline{M}_{0}
(\overline{Y},d[C_{0}])$ of genus 0 stable maps of degree $d[C_{0}]$ to
$\overline{Y}$ contains the moduli space $\overline{M}_{0} (Y,d[C_{0}])$ as
a connected component. The restriction of the virtual class
$[\overline{M}_{0} (\overline{Y},d[C_{0}])]^{vir}$ to $\overline{M}_{0}
(Y,d[C_{0}])$ defines a rational number $N_{d}$ which is the contribution
of $C_{0}$ to the degree $d[C_{0}]$ genus 0 Gromov-Witten invariant of
$\overline{Y}$.

We construct a local Calabi-Yau 3-fold $Y_{n}$ with a contractable chain of
curves $C_{n}=C_{e_{1}}\cup \dots \cup C_{e_{n}}\subset Y_{n}$ and a local
isomorphism
\[
g_{n}:Y_{n}\to Y
\]
restricting to the local immersion
\[
g_{n}:C_{n}\to C_{0}.
\]
To construct this, let $g_{\infty }:Y_{\infty }\to Y$ be the universal
cover of $Y$. We may assume that $Y$ is a sufficiently small enough
neighborhood of $C_{0}$ so that $\pi _{1} (Y)\cong \pi _{1} (C_{0})\cong
\znums $ and so $Y_{\infty }$ is an open neighborhood of $C_{\infty }$, a
linear chain of a countable number of rational curves. Fix a subchain
$C_{n}$ of length $n$ and let $Y_{n}\subset Y_{\infty }$ be an open
neighborhood of $C_{n}$ and let $g_{n}$ be the restriction of $g_{\infty }$
to $Y_{n}$.  Since $K_{Y_\infty}=g_\infty^*K_Y$, we see that $K_{Y_n}$
is trivial.

\begin{lem}\label{lem: superrigid implies Cn contracts}
Super-rigidity of $C_{0}\subset Y$ implies that $C_{n}\subset Y_{n}$ is a
generic contractible curve in the sense of subsection \ref{subsec: generic
contractable An}.
\end{lem}

\textsc{Proof:} Let $N_{C_{0}/Y}$ be the normal bundle of $C_{0}$ in $Y$
(note that $C_{0}$ and  $C_{n}$ are all local complete
intersections). By super-rigidity, 
\[
H^{0} (C_{n },g_{n }^{*}
(N_{C_{0}/Y}))=0.
\]
Since $g_{n }$ is a local isomorphism, $g_{n
}^{*} (N_{C_{0}/Y})\cong N_{C_{n }/Y_{n }}$. For any subchain
$C_{ij}\subset C_{n }$, we have a injective sheaf map 
\[
0\to N_{C_{ij}/Y_{n
}}\to N_{C_{n }/Y_{n }}|_{C_{ij}}
\]
and so $H^{0} (C_{ij},N_{C_{ij}/Y_{n
}})=0$ as well.

In particular, each component of $C_{n}\subset Y_{n}$ is a
$(-1,-1)$-curve. Thus $C_{n}$ is a contractable curve and each subchain
$C_{ij}$ has no infinitesimal deformations. Thus by Lemma \ref{lem: f
transverse is equiv to no inf defs}, we then have that $C_{n}\subset Y_{n}$
is generic.\qed

The connected components of $\overline{M}_{0} (Y,d[C_{0}])$ are identified
by the following proposition:
\begin{prop}\label{prop: components of M(C0,d)}
The connected components of $\overline{M}_{0} (Y,d[C_{0}])$ are indexed by
$n$-tuples $(d_{1},\dots ,d_{n})$ of positive integers with $\sum d_{i}=d$
and they are isomorphic to $\overline{M}_{0} (Y_{n},(d_{1},\dots ,d_{n}))$
where $(d_{1},\dots ,d_{n})\in H_{2} (Y_{n},\znums )\cong H_{2}
(C_{n},\znums )$ via the natural basis of $H_{2} (C_{n},\znums )$ indexed
by the components of $C_{n}$. Furthermore, the virtual class on
$\overline{M}_{0} (Y,d[C_{0}])$ agrees with the class induced by $Y_{n}$ on
$\overline{M}_{0} (Y_{n},(d_{1},\dots ,d_{n}))$.
\end{prop}

It follows immediately from this proposition and the genus 0 case of
Proposition \ref{prop: GW of chain} that
\[
N_{d}=\sum _{n|d}\frac{1}{n^{3}}
\]
which proves Theorem \ref{thm: formula for 1-nodal curves}.

\smallskip
\textsc{Proof:} Let $f:\Sigma \to Y$ be a genus 0 stable map. Since $\pi
_{1} (\Sigma )=0$, $f$ lifts to a map to the universal cover
\[
\til{f}:\Sigma \to Y_{\infty }.
\]
This map is unique up to deck transformations. Since $\Sigma $ is compact,
the image of $\til{f}$ is supported on some finite chain. After a deck
transformation, we can assume that the image of $\til{f}$ is exactly
$C_{n}$ for some $n$. Therefore, each $f$ factors uniquely as $g_{n}\circ
\til{f}$ where $\til{f}\in \overline{M}_{0} (Y_{n},(d_{1},\dots ,d_{n}))$
for some $n$-tuple $(d_{1},\dots ,d_{n})$ of positive integers with $\sum
d_{i}=d$. The same argument works equally well for families of stable maps
and so we get an isomorphism of moduli functors
\[
\coprod _{(d_{1},\dots ,d_{n}),\sum d_{i}=d}\overline{M}_{0}
(Y_{n},(d_{1},\dots ,d_{n})) \to \overline{M}_{0} (Y,d[C_{0}])
\]
where the isomorphism is given by composition with the appropriate
$g_{n}$. 

The virtual class $[\overline{M}_{0} (Y_{n},(d_{1},\dots ,d_{n}))]^{vir}$
is determined (as in \cite{Be-Fa}) by the perfect relative
obstruction theory
\[
[R^{\bullet}\pi _{*}\til{f}^{*} (T_{Y_{n}})]^{\vee }\to
L^{\bullet}_{\overline{M}_{0} (Y_{n},(d_{1},\dots
,d_{n}))/\mathfrak{M}_{0}}
\]
where $\pi :\mathcal{C}\to \overline{M}_{0} (Y_{n},(d_{1},\dots ,d_{n}))$
is the universal curve and $\til{f}:\mathcal{C}\to Y_{n}$ is the universal
map. Since $g_{n}$ is a local isomorphism, $T_{Y_{n}}\cong g^{*}_{n}T_{Y}$
and so $\til{f}^{*} (T_{Y_{n}})=\til{f}^{*}g_{n}^{*} (T_{Y})$. Since
$g_{n}\circ \til{f}$ is the universal map for $\overline{M}_{0}
(Y,d[C_{0}])$ under the above isomorphism, the two obstruction theories are
isomorphic and thus the virtual classes agree.\qed

\begin{rem}\label{rem: problem with higher genus}
The difficulty in generalizing this proof to higher genus Gromov-Witten
invariants is that when $\Sigma $ is no longer simply connected, we must
consider maps that factor through the finite covers of $Y$ as well. The two
terms in Conjecture \ref{conj: higher genus for 1-nodal curves} correspond
to the contributions of (1) maps factoring through the universal cover
(which we know), and (2) maps that don't factor where the conjecture
predicts the contribution is the same as if $C_{0}$ were a smooth
super-rigid elliptic curve. Note that for the genus 0 invariants, we only
need to assume 0-super-rigidity whereas to extend the argument to higher
genus we would need 1-super-rigidity as well.
\end{rem}

\begin{rem}\label{rem: problem with more nodes}
The difficulty in generalizing this proof to a rational curve with more
nodes is that the universal cover is a curve whose dual graph is a
complicated tree. One thus has to deal with curves not necessarily of ADE
type and so our deformation technique does not apply. However, for degree
five or less, the maps will factor through a map to the universal cover
that either has an ADE image or is degree one on each component. Thus under
the appropriate genericity conditions we can derive a formula for multiple
covers of $n$-nodal irreducible rational curves for multiple covers of
degree 5 or less.
\end{rem}

\section{Explicit computation in the length two ($cD_{4}$) case}\label{sec: apps and examples}

In this section, we show how to compute the multiplicities $k_{1}$ and
$k_{2}$ (see Theorem \ref{thm: multiple cover formula for embedded curves})
directly from explicit equations for the 3-fold in the case of a curve
contracting to a $cD_{4}$ singularity.

Recall that in the case of contractable a smooth rational curve $C$ inside
$Y$ , the length of the singular point $p=\pi ( C) $ inside $X$
is at most six. When the length equals one, $p$ is a $cA_{1}$ singularity
and the divisor $D^{curv}$ inside $\operatorname{PRes}( \Gamma ,\Gamma
_{0}) $ is just the origin of the affine line.

In this section we discuss in detail the length two case which corresponds
to $p$ being a $cD_{4}$ singularity. First we want to describe the divisor
$ D^{curv}=D_{1}\cup D_{2}$ in a suitable coordinate system of
$\operatorname{PRes}( \Gamma ,\Gamma _{0}) $. To start, let
$x_{1},x_{2},x_{3}$ and $x_{4}$ be any coordinate system for
$\operatorname{PRes}( \Gamma ,\Gamma _{0}) $ and the threefold $Y$ is given
by a map $f$:
\begin{eqnarray*}
f &:&\Delta \rightarrow \operatorname{PRes}( \Gamma ,\Gamma _{0})  \\
f( t)  &=&( x_{1}( t) ,x_{2}( t)
,x_{3}( t) ,x_{4}( t) ) .
\end{eqnarray*}
Notice that $x_{i}( t) $ vanishes at $t=0$ for any $i$. In fact we can
choose a good coordinate system for $\operatorname{PRes}( \Gamma ,\Gamma
_{0}) $ in such a way that $X$ can be described in term of an explicit
equation as follows (see \cite{Katz} and \cite{Ka-Mo})\footnote{ The
coordinate $x_{1}$ (resp. $x_{2},x_{3}$ and $x_{4}$) we use here
corresponds to $\sigma _{1}$ (resp. $\sigma _{2},s_{2}$ and $\sigma
_{2}^{\prime }$) in \cite{Katz}.}:
\begin{equation*}
X=\left\{ x^{2}+\frac{1}{z}( yz+x_{2}x_{4}) ^{2}=\frac{1}{z}F(
z,t) \right\} \subset \cnums ^{3}\times \Delta ,
\end{equation*}
where $F( z,t) =\left( ( z-x_{2}) ^{2}+x_{1}^{2}z \right) (
z^{2}+x_{3}z+x_{4}^{2})$ and $(x,y,z,t)$ are coordinates on $\cnums^{3}
\times \Delta $.

By direct computations, $D^{curv}=D_{1}\cup D_{2}$ can be described
explicitly as follows: 
\begin{eqnarray*}
D_{1} &=&\left\{ ( x_{2}^{2}+x_{2}x_{3}+x_{4}^{2})
^{2}+x_{1}^{2}(
x_{1}^{2}x_{4}^{2}-4x_{2}x_{4}^{2}-x_{3}x_{4}^{2}-x_{3}x_{2}^{2})
=0\right\} , \\
D_{2} &=&\left\{ x_{1}=0\right\} .
\end{eqnarray*}

From a local analysis near the singular points of the partial resolution of a
$ D_{4}$ singularity, it is shown in \cite{Katz} that the smoothness of
$Y$ implies that the order of vanishing of $x_{i}( t) $ at $t=0$ is exactly
one for $i=2,3$ or $4$. On the other hand, the order of vanishing of $
x_{1}( t) $ at $t=0$ is precisely $k_{2}$ via the above description of
$D_{2}$.

To determine $k_{1}$ we must compute the intersection number of $
f( \Delta ) $ with the hypersurface $D_{1}$ at the origin. In
fact $k_{1}$ is at least four since $D_{1}$ is defined by a polynomial which
vanishes at the origin to the fourth order. We expect that $k_{1}=4$  in
generic situation. We can describe this generic condition rather easily in
terms of either $F( z,t) $ or $f$. If we set $z=
x_{2}( t) $ then we have 
\begin{equation*}
F( x_{2}( t) ,t) =x_{1}^{2}x_{2}(
x_{2}^{2}+x_{2}x_{3}+x_{4}^{2}) 
\end{equation*}
Since $x_{1}( t) $ vanishes at $t=0$ to order $k_{2}$ and the order of
vanishing of the other $x_{i}( t) $'s is one, the order of vanishing of $F(
x_{2}( t) ,t) $ at $t=0$ is at least $2k_{2}+3.$ If it is exactly
$2k_{2}+3$, then $f( \Delta ) $ is transverse to the hypersurface $
x_{2}^{2}+x_{2}x_{3}+x_{4}^{2}=0$. By the description of $D_{1}$, this
implies that $k_{1}=4$. An alternative way to see this is as follows. Let us
write $x_{i}( t) =c_{i}t+h_{i}( t) t^{2}$ for $i=2,3$ or $4$ for nonzero
$c_{i}$'s. Then $c=\left[ c_{2},c_{3},c_{4} \right] $ determines a point in
$\mathbf{P}^{2}$ . Then we have $k_{1}=4$ as long as this point does not
lie inside the conic defined by $
x_{2}^{2}+x_{2}x_{3}+x_{4}^{2}=0$. Therefore in such generic situation, we
have $k_{1}=4$ and $k_{2}$ equals of order of vanishing of $x_{1}( t) $ at
$t=0$. When $c$ does lie inside the above conic, then, in order to
determine $k_{1}$, we also need to consider the fifth order terms in the
defining equation of $D_{1}$.  We would still obtain an explicit
description of $k_{1}$ but it will be not as simple as in the generic case.

\smallskip
\textbf{Acknowledgements:} The authors are pleased to acknowledge helpful
conversations with A.~Bertram, I.~Dolgachev, R.~Pandharipande, and C.~Vafa.

%\bibliography{mainbiblio}

\begin{thebibliography}{10}

\bibitem{Aspinwall-Morrison}
Paul~S. Aspinwall and David~R. Morrison.
\newblock Topological field theory and rational curves.
\newblock {\em Comm. Math. Phys.}, 151(2):245--262, 1993.

\bibitem{BPV}
W.~Barth, C.~Peters, and A.~Van de~Ven.
\newblock {\em Compact Complex Surfaces}, volume~4 of {\em Ergebnisse der
  Mathematik und ihrer Grenzgebiete}.
\newblock Springer-Verlag, 1984.

\bibitem{Be-Fa}
K.~Behrend and B.~Fantechi.
\newblock The intrinsic normal cone.
\newblock {\em Invent. Math.}, 128(1):45--88, 1997.

\bibitem{Br1}
E.~Brieskorn.
\newblock Singular elements of semi-simple algebraic groups.
\newblock pages 279--284, 1971.

\bibitem{Br2}
Egbert Brieskorn.
\newblock \"{U}ber die {A}ufl\"osung gewisser {S}ingularit\"aten von
  holomorphen {A}bbildungen.
\newblock {\em Math. Ann.}, 166:76--102, 1966.

\bibitem{Br3}
Egbert Brieskorn.
\newblock Die {A}ufl\"osung der rationalen {S}ingularit\"aten holomorpher
  {A}bbildungen.
\newblock {\em Math. Ann.}, 178:255--270, 1968.

\bibitem{Br-Le1}
Jim Bryan and Naichung~Conan Leung.
\newblock The enumerative geometry of {$K3$} surfaces and modular forms.
\newblock alg-geom/9711031, 1997.

\bibitem{Cl}
Herbert Clemens.
\newblock The infinitesimal {A}bel-{J}acobi mapping and moving the $\mathcal
  {O}(2)\oplus\mathcal {O}(-4)$ curve.
\newblock {\em Duke Math. J.}, 59(1):233--240, 1989.

\bibitem{CKM}
Herbert Clemens, J{\'a}nos Koll{\'a}r, and Shigefumi Mori.
\newblock Higher-dimensional complex geometry.
\newblock {\em Ast\'erisque}, (166):144 pp. (1989), 1988.

\bibitem{coxkatz}
David~A. Cox and Sheldon Katz.
\newblock {\em Mirror symmetry and algebraic geometry}.
\newblock American Mathematical Society, Providence, RI, 1999.

\bibitem{Fa-Pa}
C.~Faber and R.~Pandharipande.
\newblock Hodge integrals and {G}romov-{W}itten theory, 1998.
\newblock Preprint, math.AG/9810173.

\bibitem{Fr}
Robert Friedman.
\newblock Simultaneous resolution of threefold double points.
\newblock {\em Math. Ann.}, 274(4):671--689, 1986.

\bibitem{Go-Va}
Rajesh Gopakumar and Cumrun Vafa.
\newblock M-theory and topological strings--{II}, 1998.
\newblock Preprint, hep-th/9812127.

\bibitem{Ji}
Jes{\'u}s Jim{\'e}nez.
\newblock Contraction of nonsingular curves.
\newblock {\em Duke Math. J.}, 65(2):313--332, 1992.

\bibitem{Katz}
Sheldon Katz.
\newblock Small resolutions of {G}orenstein threefold singularities.
\newblock In {\em Algebraic geometry: Sundance 1988}, pages 61--70. Amer. Math.
  Soc., Providence, RI, 1991.

\bibitem{Ka-Mo}
Sheldon Katz and David~R. Morrison.
\newblock Gorenstein threefold singularities with small resolutions via
  invariant theory for {W}eyl groups.
\newblock {\em J. Algebraic Geom.}, 1(3):449--530, 1992.

\bibitem{Kont}
Maxim Kontsevich.
\newblock Enumeration of rational curves via torus actions.
\newblock In {\em The moduli space of curves (Texel Island, 1994)}, volume 129
  of {\em Progr. Math.}, pages 335--368. Birkh\"auser Boston, Boston, MA, 1995.

\bibitem{Li-Tian}
Jun Li and Gang Tian.
\newblock Virtual moduli cycles and {G}romov-{W}itten invariants of general
  symplectic manifolds.
\newblock In {\em Topics in symplectic $4$-manifolds (Irvine, CA, 1996)}, pages
  47--83. Internat. Press, Cambridge, MA, 1998.

\bibitem{L-L-Yau}
Bong~H. Lian, Kefeng Liu, and Shing-Tung Yau.
\newblock Mirror principle. {I}.
\newblock {\em Asian J. Math.}, 1(4):729--763, 1997.

\bibitem{man}
Yu.~I. Manin.
\newblock Generating functions in algebraic geometry and sums over trees.
\newblock In R.~Dijkgraaf, C.~Faber, and G.~van~der Geer, editors, {\em The
  moduli space of curves}, pages 401--417. Birkhauser, 1995.

\bibitem{Pa}
R.~Pandharipande.
\newblock Hodge integrals and degenerate contributions, 1998.
\newblock Preprint, math.AG/9811140.

\bibitem{Pi}
Henry~C. Pinkham.
\newblock Factorization of birational maps in dimension $3$.
\newblock In {\em Singularities, Part 2 (Arcata, Calif., 1981)}, pages
  343--371. Amer. Math. Soc., Providence, RI, 1983.

\bibitem{Pagota}
Miles Reid.
\newblock Minimal models of canonical $3$-folds.
\newblock In {\em Algebraic varieties and analytic varieties (Tokyo, 1981)},
  pages 131--180. North-Holland, Amsterdam, 1983.

\bibitem{Teodorescu-thesis}
Titus Teodorescu.
\newblock {\em Semistable torsion-free sheaves over curves of arithmetic genus
  one}.
\newblock PhD thesis, Columbia University, 1999.

\bibitem{Va}
Israel Vainsencher.
\newblock Enumeration of $n$-fold tangent hyperplanes to a surface.
\newblock {\em J. Algebraic Geom.}, 4(3):503--526, 1995.

\bibitem{Voisin}
Claire Voisin.
\newblock A mathematical proof of a formula of {A}spinwall and {M}orrison.
\newblock {\em Compositio Math.}, 104(2):135--151, 1996.

\bibitem{Wahl}
Jonathan~M. Wahl.
\newblock Equisingular deformations of normal surface singularities. {I}.
\newblock {\em Ann. of Math. (2)}, 104(2):325--356, 1976.

\bibitem{Zerger-thesis}
Thomas Zerger.
\newblock Contracting rational curves on smooth complex threefolds, 1996.
\newblock Thesis, Oklahoma State University.

\end{thebibliography}
%\bibliographystyle{plain}

\end{document}